\documentclass{amsart}
\usepackage{amsmath,amssymb,hyperref}
\newtheorem{thm}{Theorem}%[section]
\newtheorem{lem}{Lemma}%[section]

\newtheorem{prop}[lem]{Proposition}
\newtheorem{rem}{Remark}
\newenvironment{prf}{\noindent {\bf Proof.} }{\endprf\par}
\def \endprf{\hfill  {\vrule height6pt width6pt depth0pt}\medskip}

\newcommand{\ep}{\epsilon}
\newcommand{\ga}{\gamma}\newcommand{\la}{\lambda}
\newcommand{\si}{\sigma}

\newcommand{\pt}{\partial_t}

\newcommand{\dd}{\mathrm{d}}

\newcommand{\beq}{\begin{equation}}\newcommand{\ee}{\end{equation}}

\newcommand{\bea}{\begin{eqnarray}}\newcommand{\eea}{\end{eqnarray}}

\numberwithin{equation}{section}

\begin{document}

\title[ILP for SLW with low regularity]{Ill-Posedness for Semilinear
Wave Equations with Very Low Regularity}

\author{Chengbo Wang}
\address{Department of Mathematics, Zhejiang University, Hangzhou,
310027, China} \email{wangcbo@yahoo.com.cn}
\thanks{The authors were partially supported by NSF of China 10571158.}

\author{ Daoyuan Fang}
\address{Department of Mathematics, Zhejiang University, Hangzhou,
310027, China} \email{dyf@zju.edu.cn}
%\thanks {The second author was supported by NSF of China 10571158.}

\date{}
\dedicatory{} \commby{}

\begin{abstract}
In this paper, we study the ill-posdness of the Cauchy problem for
semilinear wave equation with  very low regularity, where the
nonlinear term depends on $u$ and $\partial_t u$. We prove a
ill-posedness  result for the ``defocusing" case, and give an
alternative proof for the supercritical ``focusing" case, which
improves our previous result in Chin. Ann. Math. Ser. B 26(3),
361--378, 2005.
\end{abstract}

\maketitle

\section{Introduction}\label{intro}

This paper is mainly concerned with a low regular ill-posedness(ILP)
of the Cauchy problem for the ``defocusing" semilinear wave equation
with nonlinear term depending on $u$ and $\partial_t u$. For the
``focusing" case, it has been studied in \cite{FngWng1}, but here we
give an alternative proof of the supercritical ill-posed result in
the final section, which can improve the original result slightly.

Recently, the study of ill posed issues for nonlinear evolution
equations is very active. For wave equations, one can refer to
\cite{brenner}-\cite{95kuksin}. In \cite{brenner}, \cite{ChCoTa03},
\cite{lebeau1}, \cite{lebeau2}, \cite{95LdSo} and \cite{95kuksin},
these authors give the ill-posed results for both focusing and
defocusing type equations with nonlinear term depending only on $u$
itself.

For nonlinear term depends also on the derivatives of $u$, only a
few results can be found but for the wave map type equations, which
one can refer to \cite{DaGe05} and references therein. In
\cite{93Lind} and \cite{FngWng1}, the authors deal with the
``focusing" type semilinear equations, where ``focusing" means that
the corresponding time-ODE equation admits finite time blow up
solution. In \cite{96Ld} and \cite{Lind98}, the author dealt with
some particular``focusing"  semilinear and quasilinear equations.

Let $\Box :=\pt^2 -\Delta_x ,$ $x\in \mathbb{R}^n,$ be the usual
d'Alembertian, we consider the Cauchy problem for the following
``defocusing'' equations($k+l>1$, $k\in\mathbb{N}$ and $1\leq
l\in\mathbb{R}$)
\begin{equation}\label{ilp}\Box u=-|u|^k |\partial_t u|^{l-1}
\partial_t u,\,\, \mbox{for even}\,k\end{equation}
and \begin{equation}\label{ilp'}\Box u=- |u|^{k-1}u |\partial_t
u|^l,\,\, \mbox{for odd} \,k.
\end{equation}
If we denote the right hand side of \eqref{ilp} or \eqref{ilp'} by
$F(u,\pt u),$ then the corresponding ordinary differential
equation(ODE) in $t$ is
\begin{equation}\label{ode}
\pt^2 u=F(u,\pt u)\ .
\end{equation}

Heuristically, for ``focusing'' equations of these type, there are
two obstructions for well-posedness, one is scaling, which yields
the scaling index $s_c(k,l)=\frac{n}{2}+\frac{l-2}{k+l-1}$, another
is concentration along light cone, which yields the concentrative
index $\tilde{s}(k,l)=\frac{n+1}{4}+\frac{l-1}{k+l-1}$. For
``defocusing'' equations, the mechanism for ill posedness are
currently not very clear.

In contrast, for both ``focusing'' and ``defocusing'' equations, if
$k,l\in\mathbb{N}$ and $k+l>1$, we have local well-posedness (LWP)
in $C_t H^s$ for $$s>\left\{\begin{array}{ll}\max(\tilde{s}(0,l),
s_c(0,l)),& l \geq 2\\%\ and\ n\geq 2\\
\max(\frac{n-1}{2}, s_c(k,1)),%\vee\frac{n+1}{4}
&l=1\ \textrm{and}\ n\geq 3\end{array}\right.$$ In particular, if
$l-1\ge \frac{4}{n-1}$ and $s>s_c(0,l)$, or $k\ge 2$, $n\ge 3$ and
$s>s_c(k,1)$, we have LWP in $C_t H^s$ (see \cite{FngWng1}).
Moreover, for $k=0$ and $l-1> \max(2,\frac{4}{n-1})$, we have global
well-posedness in $C_t H^s$ with small data for $s\ge s_c(0,l)$ (see
\cite{FngWng2}).

In this paper we will mainly prove the following ill posedness
result:
\begin{thm}\label{result}
If $k+l$ is odd or $l\geq \max(k+1,[\frac{n+4}{2}])$, the problem
\eqref{ilp} or \eqref{ilp'} is ill posed in $H^s$ for
$s\in(-\infty,\frac{1}{2}]\cap
(-\infty,\frac{1}{2}+\frac{l-2}{k+l-1})$, in the sense of that the
solution map is not continuous from $H^s\times H^{s-1}$ to $C_t
H^s\times C^1_t H^{s-1}$. Precisely, we get that for any $\epsilon$,
there exists a solution $u\in C_t H^s\times C^1_t H^{s-1}$ and some
$t\in(0,\epsilon)$ such that
\begin{equation}
u(0)=0,\ \partial_t u(0)=\varphi\in\mathcal{C}^\infty_0,\
\|\varphi\|_{H^{s-1}}<\epsilon \ \mathrm{and}\ \|\partial_t
u(t)\|_{H^{s-1}}>\frac{1}{\epsilon}\ .
\end{equation}
\end{thm}

Inspired by the paper \cite{ChCoTa03} of Christ, Colliander and Tao,
we first use small dispersion analysis and scaling argument to get a
well-controlled two-parameter solution. By choosing the parameters
properly, one can get the desired estimate. Thus we have the
following result,
\begin{prop}\label{weak}
The problem \eqref{ilp} or \eqref{ilp'} is ill posed in $H^s$ for
$s\in(-\infty,\frac{2-n}{2}]\cap
(-\infty,\frac{n}{2}+\frac{l-2}{k+l-1})$ if $k+l$ is odd or $l\geq
k+1$.
\end{prop}
Then, by the argument of finite speed of propagation, Theorem
\ref{result} can be reduced to the $n=1$ case of Proposition
\ref{weak}.

\begin{rem} Our
original purpose is to get the ill posedness result for any $s<s_c$.
%, where $s_c$ is the critical regularity.
However, because of the lack
of the property of blow up at infinity for the ODE solution, we
couldn't get ill posed result by exploring the continuous
dependence. But we believe that there must be some other mechanism
to develop ill posed result.
\end{rem}

As a complementarity, for the supercritical ill-posed result of
``focusing" equation in \cite{FngWng1}, we give an alternative
proof, which can improve the original result slightly. Combined with
the result of Theorem 1.2 in \cite{FngWng1}, we'll prove the
following result in section \ref{appendix}.
\begin{thm}\label{foc}
Let $n\geq 1$, $k+l>1$, $k,l\geq 0$ and $k,l\in\mathbb{R}$. Consider
the model equation \beq\Box u=|u|^k |\partial_t u|^{l-1}\partial_t u
\ \ \ \mbox{if}\ \ l\neq 0\ee \beq \Box u=|u|^{k-1}u \ \ \
\mbox{if}\ \ l=0\ee they are s-ILP in $\Dot{H}^s$ for \beq
s\in\left\{\begin{array}{ll}(1-n/2,s_c)&l\ge
2\\(-n/2,s_c)&l<2\end{array}\right.\ee And, if $s_c> 1$, then it is
s-ILP in $H^s$ for any $s<s_c$. Moreover, it's also w-ILP in $H^s$
for $s< s_c$ if $s_c\geq 1$.
\end{thm}
Here s-ILP means that there is a sequence of data $f_j, g_j\in
C^\infty_0(B_{R_j})$, for which the lifespan of the solutions $u_j$
tends to zero as the data's norm and $R_j$ goes to $0$, under the
condition that the solutions obey finite speed of propagation; and
w-ILP to express that the lifespan goes to zero and the data's norm
stay bounded. In fact, in the case of s-ILP, one also can find that
the solution does not depends continuously on the data.

The main differences for Theorem \ref{foc} from  Theorem 1.2 in
\cite{FngWng1} are in two cases, the first is that there isn't the
technical restriction $k=0$ for $l=2$, the second is that we have
ILP in $H^s$ even for $s<-\frac{n}{2}$.

In this paper, we will use the following conventions. We use the
notation $(f,g)$ to stand for the specification of the data
$u(0)=f$, $\partial_t u(0)=g$. Let $\alpha=\frac{l-2}{k+l-1}$,
$s_c=\frac{n}{2}+\alpha$ and $B_r = \{ |x|<r \}$. Moreover, we use
$\hat{f}=\mathcal{F}(f)$ to denote the Fourier transform of the
function $f$.

\section{Proof of Proposition \ref{weak}}\label{WeakProp}
In this section, we use the so-called ``small dispersion analysis"
in \cite{ChCoTa03} to prove the ill posed result Proposition
\ref{weak}, based on the knowledge of the properties of the
solutions of ODE \eqref{ode}.

\subsection{ODE Solution}\label{ODE}
In this subsection, we study the asymptotic properties of the
solutions of ODE \eqref{ode}, which is the basis of the small
dispersion analysis.

Note that there is  a ``conserved quantity" for ODE
 \eqref{ode} (for $u,\partial_t u\geq 0$),
\begin{equation}\label{ConQuantity}
  I=\left\{
  \begin{array}{ll}
  \frac{|\partial_t
  u|^{2-l}}{2-l}+\frac{|u|^{k+1}}{k+1},& l\neq 2,\\
 \log(|\partial_t u |)+\frac{|u|^{k+1}}{k+1},& l=
 2, \,\partial_t u > 0. \end{array}\right.
\end{equation}
Combine it with the equation, it is easy to get the following global
asymptotic properties of the solution $u_1(t)$ of \eqref{ode} with
data $(0,1)$.

\begin{prop}\label{odeprop}Let $l\geq 1$, $k\geq 0$, $k+l>1$ and $M$ be
the prescribed large numbers.
  The solution $u_1(t)$ of  \eqref{ode}
  with data $u(0)=0$, $\partial_t u(0)=1$ has the following
  properties:
  \begin{enumerate}
    \item $u_1(t)$ exists globally;
    \item $\lim_{t\rightarrow\infty}u_1(t)=\left\{
        \begin{array}{ll}
        \infty& l\geq 2,\\
         \left(\frac{k+1}{2-l}\right)^{\frac{1}{k+1}}& l< 2\end{array}\right.;$
    \item
$u_1(t)\in C^{m_0+2}$ with $m_0=\left\{\begin{array}{ll}M &
k,l\in\mathbb{Z}, \ k+l\ \rm{odd},\\\min([k-1],[l-1])&
\rm{else}\end{array}\right.;$
    \item $\lim_{t\rightarrow\infty} \pt^m u_1(t)=0,\ \forall\ 0<m\leq
    m_0+2.$
  \end{enumerate}
\end{prop}

\begin{prf}
Since we have the ``conserved quantity" \eqref{ConQuantity} for the
solution, there exists a $T>0$ such that
$$\pt u_1\searrow 0,\ u_1\nearrow a\ :=\ \left\{\begin{array}{lr}
  \infty, & l\geq 2, \\
  \left(\frac{k+1}{2-l}\right)^{\frac{1}{k+1}}, & l<2
\end{array}\right.
$$ as $t\rightarrow T.$ Moreover,  we claim that for $l\geq 1$, we have $T=\infty$.
In fact, by \eqref{ConQuantity}, we get that for $t\in[0,T)$,
$$\pt u_1=f(u_1):=\left\{\begin{array}{lr}
  \left(1+\frac{l-2}{k+1}u_1^{k+1}\right)^{\frac{1}{2-l}}, & l\neq 2, \\
  \exp\left(-\frac{u_1^{k+1}}{k+1}\right), & l=2.
\end{array}\right.$$
Then, note that $f(u_1)\leq 1$ for $l\geq 2$, we have $\pt u_1\leq
1$. Thus $u_1(t)\leq t$, for any $t$, implies $T=\infty$. For $1\leq
l<2$, note that $0\leq a^{k+1}-u_1^{k+1}\leq (k+1)a^{k} (a-u_1)$ and
$\frac{1}{l-2}\leq -1$, we have
$$T=\int_0^T \dd t =c\int_0^a
\left(a^{k+1}-u_1^{k+1}\right)^{\frac{1}{l-2}}\dd u_1\geq c \int_0^a
\left(a-u_1\right)^{\frac{1}{l-2}}\dd u_1=\infty\ .
$$

From the fact that $F(a,b)\in C^{m_0,1}([0,u_1(T)]\times (0,1])$ for
any $T\in (0,\infty)$, we have $u_1(t)\in C^{m_0+2}([0,\infty))$.

At last, to prove the property of $u_1^{(m)}:=\pt^m u_1$, we divide
it into three cases $l<2$, $l>2$ and $l=2$. For $l<2$, it is obvious
since $u_1^{(m)}=P((u_1^{(j)})_{j<m}))$, where $P$ is polynomial and
$P|_{(u_1^{(j)})_{0<j<m}=0}=0$. For $l>2$ and $t$ large, $$\pt
u_1=\left(1+\frac{l-2}{k+1}u_1^{k+1}\right)^{\frac{1}{2-l}}\simeq
u_1^{\frac{k+1}{2-l}} \Rightarrow u_1\simeq t^{\frac{l-2}{k+l-1}}\
,$$ thus we have, for any integer $m>0$, $u_1^{(m)}\simeq
t^{\frac{l-2}{k+l-1}-m}\rightarrow 0$ as $t\rightarrow \infty$. For
$l=2$ and the large $t$, we have $$u_1(t)\geq \epsilon (\ln
t)^{\frac{1}{k+1}}$$ and $u_1^{(m)}(t)\leq t^{-\delta}$ for some
$\delta>0$. This completes the proof. \end{prf}

For any $\phi>0$, we can get the solution $\tilde{u}_\phi(t)$ of
\eqref{ode} with data $(0,\phi)$ by rescaling
\begin{equation}
\tilde{u}_\phi(t)=u_1(t\phi^{\frac{k+l-1}{k+1}})\phi^{-\frac{l-2}{k+1}}.
\end{equation}
For the continuity in parameter, we choose $N$ large and let
$\phi=\psi^{N(k+1)}$, \beq
u_\psi(t)=\tilde{u}_\phi(t)=u_1(t\psi^{N(k+l-1)})\psi^{-N(l-2)}\
.\ee Note that $\lim_{\psi\rightarrow 0+}u_\psi(t)=0$, we use the
convention that $u_0=0$. Then we claim that \beq\label{ContOnPara}
u_{\psi}(t)\in C_{t,\psi}^{m_0+2}([0,\infty)\times[0,\infty))\ \rm
if\ N(k+1)>m_0+2\ ,\ee and consequently $u_{\psi(x)}(t)\in
C_{t,x}^{m_0+2}$ for $\psi(x)\in C^\infty_0(\mathbb{R}^n)$ and
$\psi\geq 0$.

Now we give the proof of the claim. For the simplicity of notation,
we denote here that $f=u_1$, $b=N(k+l-1)$, $a=N(l-2)$ and
$u(t,\psi)=f(t \psi^b)\psi^{-a}$. We extend the definition of
$u(t,\psi)$ from $\psi\geq 0$ to $\psi\in\mathbb{R}$ by using the
zero extension at first. Then we check that such extension lies in
$C_{t,\psi}^{m_0+2}$ and show that the only case to be examined is
the case $\psi=0$. We finish the proof by computing the right limit
at $\psi=0$. Note that if $t>0$ and $b-d>0$, $$\lim_{\psi\rightarrow
0+}f(t \psi^b)\psi^{-d}=0\ .$$ Since we have
$$\partial_{\psi}^j u(t,\psi) =  C f(t \psi^b)\psi^{-a-j}+
\sum_{1\le h\le i\le j} C_{h,i} t^h (\partial^h f)(t \psi^b) \psi^{b
h -a-j}$$ then if $b-a-j>0$, we have
\begin{equation}\label{limit}
\lim_{\psi\rightarrow 0+}\partial_{\psi}^j u(t,\psi)=0\
.\end{equation} Thus we have \eqref{limit} for any $0\le j\le m_0+2$
if $N(k+1)=b-a>m_0+2$. This complete the proof of the claim.

\subsection{Small Dispersion Analysis}\label{SDA}

Based on Proposition \ref{odeprop} and \eqref{ContOnPara}, we make
the small dispersion analysis in this subsection.

Consider the problem \beq\label{sda} \left\{
\begin{array}{l}\Box_\ga u :=(\partial_t^2-\ga^2\Delta)
u=F(u,\pt u)\\
u(0)=0,\ \pt u(0)=\phi(x)=\psi(x)^{N(k+1)}
\end{array}\right.
\ee By the usual energy argument, we can compare the solution
$\phi^{(\gamma)}$ with the corresponding solution
$\phi^{(0)}=u_\psi$ of \eqref{ode}. The result is as following
\begin{prop}\label{prop}
Let $n\geq 1$, $k \geq 0$, $l \geq 1$, $k+l > 1$, $m\in \mathbb{N}$,
$[\frac{n+2}{2}]\leq m\leq m_0$, $N$ such that $N(k+1)>m+2$,
$\phi=\psi^{N(k+1)}$ with $\psi\in C^\infty_0(\mathbb{R}^n)$ and
$\psi\geq 0$. Then there exist $C>>1>>c>0$ such that $\forall
\gamma\in(0,c]$, there exist a solution $u=\phi^{(\gamma)}\in
C([0,T],H^{m+1})\cap C^1([0,T],H^{m})$ of \eqref{sda} with
$T=c|\log(\gamma)|^c$. Moreover, for any $t\in[0,T]$,
\begin{equation}
\|\phi^{(\gamma)}-\phi^{(0)}\|_{H^{m+1}}+\|\partial_t
\phi^{(\gamma)}-\partial_t \phi^{(0)}\|_{H^{m}}\leq C\gamma^{1/2}.
\end{equation}
\end{prop}

\begin{prf}
  Let $w=u-\phi^{(0)}$, then we have
\begin{equation}\label{ep3}
\Box_\gamma w=F(u,\partial_t u)-F(\phi^{(0)},\partial_t
\phi^{(0)})+\gamma^2 \Delta \phi^{(0)}=G(w)
\end{equation} for $w$ with data
$(0,0)$.

The energy method shows that this problem is local well-posed in
$H^{m+1} \times H^{m}$, the solution of \eqref{ep3} exists as long
as the $H^{m+1} \times H^m$ norm of it stays bounded.

We define the $\gamma$-energy of $w$ by
$$ E_\gamma(w(t)) := \int \frac{1}{2} |w_t(t,y)|^2 + \frac{\gamma^2}{2}
|\nabla_y w(t,y)|^2\ dy\ ,$$
$$ E_{\gamma,m}(w(t)) := \sum_{j=0}^m E_\gamma(\partial_y^j w(t))\ .$$
Then the standard energy inequality gives that
\begin{equation}\label{energy} |\partial_t
E_{\gamma,m}^{1/2}(w(t))| \leq C \| G(t) \|_{H^m}.
\end{equation}

Let $e(t)$ be the non-decreasing function $ e(t) := \sup_{0 \leq s
\leq t} E_{\gamma,m}^{1/2}(w(s))$, then
 \begin{equation}\label{w-bound} \| w(t) \|_{H^m}
\leq \int_0^t \| w_t(s)\|_{H^m}\ ds \leq C \int_0^t
E_{\gamma,m}^{1/2}(w(s))\ ds \leq C t e(t),
\end{equation}
Also, by the smoothness of $\phi^{(0)}$ and $F$, we can easily
obtain the bounds
$$ \| \gamma^2 \Delta_y \phi^{(0)} \|_{H^m} \leq C \gamma^2 (1 + |t|)$$
and
$$ \| \phi^{(0)} \|_{H^m} + \| \partial_t \phi^{(0)} \|_{H^m} + \| \phi^{(0)}
\|_{C^m} + \|\partial_t \phi^{(0)} \|_{C^m} \leq C (1 + |t|)\ .$$
Since $H^m$ is an algebra, then
$$ \| F(u, \partial_t u)(t) - F(\phi^{(0)},\pt \phi^{(0)})(t)\|_{H^m}
\leq C (1+|t|)^C(e(t)+e(t)^C)\ .$$

The above estimates yield that
$$ \| G \|_{H^m}
\leq C (1 + |t|)^C ( \gamma^2 + e(s) + e(s)^C ),$$ which by
\eqref{energy} gives the differential inequality
$$ \partial_t e(t) \leq C (1 + |t|)^C (\gamma^2 + e(t) + e(t)^C).$$

Since $e(0) = 0$, we can assume {\it a priori} that $e(t)\leq\ga$,
then
$$ \partial_t e(t) \leq C (1 + |t|)^C (\gamma^2 + e(t)),$$
and hence $$e(t)\leq \gamma^2 (\exp(C(1 + |t|)^C)-1)\ .
$$
Thus if $|t| \leq c |\log \gamma|^c$ for suitably chosen $c$ and
$\gamma$ to be sufficiently small, we obtain
$$ e(t) \leq C \gamma^{3/2}$$
and furthermore we can recover the {\it a priori} assumption, which
can then be removed by the usual continuity argument. The claim then
follows from \eqref{w-bound} if $\gamma$ is sufficiently small.
\end{prf}

\subsection{Estimate for Solution}\label{sltest}
By Proposition \ref{prop} and rescaling ($\la>0$), we get
two-parameter solutions for the problem \eqref{sda} with $\gamma=1$,
%(i.e. the equation \eqref{ilp} or \eqref{ilp'})
\begin{equation}
u^{(\gamma,\lambda)}:=\lambda^\alpha
\phi^{(\gamma)}(\lambda^{-1}t,\lambda^{-1}\gamma x)\ .
\end{equation}
In particular, we have the initial data
$$(0,\lambda^{-\frac{k+1}{k+l-1}}\phi(\lambda^{-1}\gamma x))\ .$$

\begin{prop}\label{prop'}
Let $0<\lambda\leq\gamma\ll 1$, \begin{equation}   \|\partial_t
u^{(\gamma,\lambda)}(0)\|_{H^{s-1}}=
C\lambda^{s_c-s}\gamma^{s-\frac{n+2}{2}}:=C\epsilon
\end{equation}
where $\phi \in C^\infty_0(\mathbb{R}^n)$ such that
$\hat{\phi}(\xi)=O(|\xi|^k)$ as $\xi\rightarrow 0$ with $k+s>1-n/2$.
Note that $\la=c\ga^{\si}$ with $\si>1$ for fixed
$\ep$ and $s<s_c$.% if $s\leq 1-n/2$.
\end{prop}

\begin{prf}

Note that
  $$  [\partial_t u^{(\gamma,\lambda)}(0)]\widehat{\phantom{A}}(\xi) =
\lambda^{-\frac{k+1}{k+l-1}} ({\lambda}/{\gamma})^n
\widehat{\phi}({\lambda \xi}/{\gamma})\ ,$$ we have
\begin{align*}
\|\partial_t u^{(\gamma,\lambda)}(0)\|_{H^{s-1}}^2 &=
\lambda^{-2\frac{k+1}{k+l-1}} (\lambda/\gamma)^{2 n} \int
|\widehat{\phi}(\lambda\gamma^{-1}\xi)|^2 (1+|\xi|^2)^{s-1}\,d\xi
\\
&= \lambda^{-2\frac{k+1}{k+l-1}} (\lambda/\gamma)^{n} \int
|\widehat{\phi}(\eta)|^2
(1+|\gamma\lambda{^{-1}}\eta|^2)^{s-1}\,d\eta.
\\
&\sim \lambda^{-2\frac{k+1}{k+l-1}} (\lambda/\gamma)^{n-2(s-1)}
\int_{|\eta|\ge\lambda\gamma{^{-1}}}
|\widehat{w}(\eta)|^2|\eta|^{2(s-1)}\,d\eta
\\
& \qquad\qquad+ \lambda^{-2\frac{k+1}{k+l-1}} (\lambda/\gamma)^{n}
\int_{|\eta|\le\lambda\gamma{^{-1}}} |\widehat{\phi}(\eta)|^2\,d\eta
\\
&= \lambda^{-2\frac{k+1}{k+l-1}} (\lambda/\gamma)^{n-2(s-1)}\left[
\int_{\mathbb{R}^n}
|\widehat{\phi}(\eta)|^2|\eta|^{2(s-1)}\,d\eta\right.
\\
&\qquad \left.- \int_{|\eta|\leq\lambda\gamma^{-1}}
|\widehat{\phi}(\eta)|^2\,
\big((\lambda/\gamma)^{2(s-1)}-|\eta|^{2(s-1)}\big)\,d\eta\right].
\end{align*}

Thus for any $s-1>-n/2$,
\begin{equation}
\|\partial_t u^{(\gamma,\lambda)}(0)\|_{H^{s-1}} =
c\lambda^{-\frac{k+1}{k+l-1}}({\lambda}/{\gamma})^{n/2 - (s-1)}
\cdot(1 + O\big((\lambda\gamma^{-1})^{s-1+n/2}\big)),
\end{equation}
where $c\ne 0$ provided that $\phi$ is not identically zero. In
particular,
\begin{equation} \label{time0normal}
\|\partial_t u^{(\gamma,\lambda)}(0)\|_{H^{s-1}} =
c\lambda^{-\frac{k+1}{k+l-1}}({\lambda}/{\gamma})^{n/2 -
(s-1)}\end{equation} provided that $s-1>-n/2$ and
$\lambda\ll\gamma$.

For $s-1\leq -n/2$, \eqref{time0normal} still holds, under the
supplementary hypothesis that
\begin{equation}  \label{momentcondition}
\widehat{\phi}(\xi) =O(|\xi|^k) \ \text{as $\xi\to 0$, for some $k>
-(s-1)-n/2$.}
\end{equation}
Then, if $\lambda\ll\gamma$, we have $\int_{\mathbb{R}^n}
|\widehat{\phi}(\eta)|^2|\eta|^{2(s-1)}\,d\eta<\infty$ and
\begin{equation*}
\int_{|\eta|\le\lambda\gamma^{-1}} |\widehat{\phi}(\eta)|^2\,
\big((\lambda/\gamma)^{2(s-1)}-|\eta|^{2(s-1)}\big)\,d\eta \le
C(\lambda\gamma^{-1})^{n+2(s-1)+2 k} \le C<\infty.
\end{equation*}
\end{prf}

To complete the proof of Proposition \ref{weak}, we need to choose
the data $\varphi$ appropriately.

At first glance, for $s\le 1-n/2$, the condition of the data in
Proposition \ref{prop'} couldn't be fulfilled since we assume the
data to be nonnegative in Proposition \ref{prop}. However, noting
that $-u$ is also a solution of \eqref{ilp} or \eqref{ilp'} whenever
$u$ is, and the solution exhibit uniformly finite speed of
propagation so long as $|\gamma|\le 1$, such condition can be easily
fulfilled by taking $\varphi$ to be an appropriate linear
combination of nonnegative $C^\infty_0$ functions with widely spaced
supports as in \cite{ChCoTa03}.

Moreover, by the conditions about $k$ and $l$ in Proposition
\ref{weak}, we can  choose $\varphi$ such that $\varphi^{(0)}$ has
the following property:  there is a $t_0> 0$ such that
\begin{equation}
\label{con}\mathcal{F}(\partial_t \varphi^{(0)}(t_0))(0)\neq
0\end{equation} even if $$\mathcal{F}(\partial_t
\varphi^{(0)}(0))(0)= 0\ .$$

In fact, we consider, for the case that $k$ is even for example, the
quantity $A(t)=\int\partial_t \varphi^{(0)}(t,x)\mathrm{d} x$.
 Then from the equation \eqref{ode},
we have $$\pt A(t)=-\int |\varphi^{(0)}|^k |\pt
\varphi^{(0)}|^{l-1}\pt \varphi^{(0)} \mathrm{d}x\ ,$$ which is
vanishing at $t=0$ for $k>0$ since $\varphi^{(0)}$ is zero at $t=0$.
Similarly, for any $i\leq k$, $\pt^i \varphi^{(0)}$ is also null at
$t=0$. Since $k$ is even, we have
\begin{equation}\label{cond} \pt^{k+1} A(0)=k!\int
|\pt\varphi^{(0)}(0)|^{k+l-1}\pt\varphi^{(0)}(0) \mathrm{d}x=k!\int
|\varphi|^{k+l-1}\varphi \mathrm{d}x
\end{equation} at $t=0$. Then we may require that $\varphi$ satisfies an
additional condition that the right hand side of \eqref{cond} is
nonzero. Hence we get \eqref{con} immediately. Here, if $l$ isn't
odd, we use the condition $l\geq k+1$ to ensure $\varphi^{(0)}\in
C^{k+2}$.

Then, from \eqref{con}, we get for some $c>0$
$$|\mathcal{F}(\partial_t \varphi^{(0)}(t_0))(\xi)|\geq c\ \ \rm{when}\ \ |\xi|
\leq c\ .$$ However, by Proposition \ref{prop} and Sobolev
embedding,
\begin{eqnarray*}|\mathcal{F}(\partial_t
(\varphi^{(\gamma)}-\varphi^{(0)}))(\xi)|&\lesssim&
\|\mathcal{F}(\partial_t
(\varphi^{(\gamma)}-\varphi^{(0)}))(\xi)\|_{H^m}\lesssim
 \sum_{j = 0}^m \|
{|x|}^{j} \pt(\varphi^{(\gamma)}-\varphi^{(0)}) \|_{L^2}
\\&\lesssim&
\|\pt(\varphi^{(\gamma)}-\varphi^{(0)})\|_{L^2}\lesssim C\gamma\ ,
\end{eqnarray*}
where we have used the fact that
$\pt(\varphi^{(\gamma)}-\varphi^{(0)})$ are supported in a fixed
compact set.
 Thus we get that
$$|\mathcal{F}(\partial_t \varphi^{(\gamma)}(t_0))(\xi)|\geq c\ \ \rm{for}\ \ |\xi|
\leq c$$
\begin{equation}
  |\mathcal{F}(\partial_t u^{(\gamma,\lambda)}(\lambda t_0))(\xi)|\geq c
\lambda^{\alpha-1}  \left(\frac{\gamma}{\lambda}\right)^{-n}\ \
\rm{for}\ \ |\xi|\leq c \frac{\gamma}{\lambda}
\end{equation}

Note that $\gamma\gg \lambda$ for $\gamma$ small and $\epsilon$
fixed. For $s<1-n/2$ and $s<s_c$,
\begin{equation}
  \|\partial_t u^{(\gamma,\lambda)}(\lambda t_0)\|_{H^{s-1}}\geq c \lambda^{\alpha-1}
  \left(\frac{\gamma}{\lambda}\right)^{-n}= c\epsilon
  \left(\frac{\gamma}{\lambda}\right)^{1-n/2-s}=1/\epsilon
\end{equation}
if $\gamma$ is sufficiently small. For $s=1-n/2<s_c$,
\begin{eqnarray*}
  \|\partial_t u^{(\gamma,\lambda)}(\lambda t_0)\|_{H^{-n/2}}^2&
  \geq &
  c \lambda^{2\alpha-2}
  \left(\frac{\gamma}{\lambda}\right)^{-2n}\int_{|\xi|\leq
  c\frac{\gamma}{\lambda}}(1+|\xi|)^{-n}\mathrm{d}\xi\\
  &\geq& c \lambda^{2\alpha-2}
  \left(\frac{\gamma}{\lambda}\right)^{-2n}
  \log(c\frac{\gamma}{\lambda})\\
  &=&c\epsilon^2\log(c\frac{\gamma}{\lambda})=\epsilon^{-2}
\end{eqnarray*}
This complete the proof of Proposition \ref{weak}.

\section{Reduction to dimension $n=1$}\label{Red1-D}

By the argument of finite speed of propagation, Theorem \ref{result}
can be reduced to the $n=1$ case of Proposition \ref{weak}.

This reduction is worked by considering initial data of product form
$\eta(x')\varphi(x_n)$ where
$x=(x',x_n)\in\mathbb{R}^{n-1}\times\mathbb{R}$, $\varphi$ is the
same as that in the proof of Proposition \ref{weak}, and $\eta$ is a
fixed $C_0^\infty$ function which equals $1$ on a sufficiently large
ball in $\mathbb{R}^{n-1}$. By finite speed of propagation, the
corresponding solutions, assuming existence and uniqueness, will
likewise have product form for $x'$ in a fixed smaller ball, so the
norm estimate in $\mathbb{R}^n$ follows from the estimate already
established in $\mathbb{R}^1$.

Precisely, we set $\eta(x')\in C^{\infty}_0$ and $\eta=1$ on the
ball $B_R(\mathbb{R}^{n-1})$ with radius $R\gg 1$. Then for some
$\nu\ll 1$, the problem \eqref{ilp} or \eqref{ilp'} with data
$(0,\eta(x')\varphi(x_n))$ are local well-posed in
$C([0,\nu],H^{m+1})\cap C^1([0,\nu],H^{m})$ with
$m=[\frac{n+2}{2}]$(Note that this is where we need the extra
condition that $l\geq [\frac{n+4}{2}]$). Thus the solution $u(t)$
has the property of ``finite speed of propagation", and hence for
$(x',x_n,t)\in B_{R-1}(\mathbb{R}^{n-1})\times\mathbb{R}\times
[0,\nu]$, the value of $u(x,t)$  depends only on the data in
$B_{R}(\mathbb{R}^{n-1})\times \mathbb{R}$, i.e., only on
$\varphi(x_n)$. Thus we have that
$$u(x,t)=\tilde{u}_{\varphi}(x_n,t) \ \rm in\ B_{R-1}(\mathbb{R}^{n-1})
\times\mathbb{R}\times[0,\nu], $$ where $\tilde{u}_{\varphi}$ denote
the solution of \eqref{ilp} or \eqref{ilp'} with data $(0,\varphi)$.

Now we give the corresponding estimate of such data and solution.
Note that $s-1<-n/2<0$ and $|\xi|\geq|\xi_n|$, we have
$$\|\eta(x')\varphi(x_n)\|_{H_{x}^{s-1}}\leq
C\|\varphi\|_{H_{x_n}^{s-1}}\|\eta\|_{L_{x'}^2}\leq C\epsilon\ .$$
Let $f(x)=g(x')h(x_n)\in C^{\infty}_0(\mathbb{R}^{n)}$ with $h(x_n)
\tilde{u}_{\varphi}(x_n,t)=\tilde{u}_{\varphi}(x_n,t)$, $g$
supported in $B_{R-1}(\mathbb{R}^{n-1})$ and
$|\widehat{g}(\xi')|\geq c>0$ for $\xi' \le c$. Since the
generalized Leibniz rule(see e.g. Lemma 2.2 of \cite{FngWng1})
yields that for any $s\in\mathbb{R}$ and $\epsilon>0$,
\begin{equation}\label{Leib}
\|f u\|_{H^s}\leq C \|f\|_{H^{\max(|s|,
\frac{n}{2}+\epsilon)}}\|u\|_{H^s}=C_f \|u\|_{H^s}\ ,\end{equation}
we have
\begin{eqnarray*}
  \|\pt u(t_0)\|_{H^{s-1}}&\geq&\frac{1}{C_{f}}\|f \pt u(t_0)\|_{H^{s-1}}\\
&=&  \frac{1}{C_{f}}\|g(x') h(x_n) \pt \tilde{u}_{\varphi}(x_n,t_0)\|_{H_x^{s-1}}\\
  &\geq& c\| h(x_n) \pt \tilde{u}_{\varphi}(x_n,t_0)\|_{H_{x_n}^{s-1}}\\
%Here $\geq$ since we can limit the integrate domain to $|\xi'|\leq c$.
  &=&c\|\pt \tilde{u}_{\varphi}(x_n,t_0)\|_{H_{x_n}^{s-1}} \ \geq \ c \epsilon^{-1}
\end{eqnarray*}
This complete the reduction.

\section{An alternative proof of $s<s_c$ ILP in
\cite{FngWng1}}\label{appendix}

In this section, combined with Theorem 1.2 in \cite{FngWng1}, we
prove Theorem \ref{foc} for the following ``focusing" equation:
\beq\label{focusing-eqn}\Box u=|u|^k |\partial_t u|^{l-1}\partial_t
u \ \ \ \mbox{if}\ \ l\neq 0\ee \beq \Box u=|u|^{k-1}u \ \ \
\mbox{if}\ \ l=0\ee with $k+l>1$, $k,l\geq 0$ and
$k,l\in\mathbb{R}$.

Several supercritical ill posed(ILP) results of above  equations
have been obtained in our previous paper \cite{FngWng1}, here we
give an alternative proof  and a slightly improvement for Theorem
1.2 in there. In \cite{FngWng1}, the starting point is the explicit
blow-up solution in time-ODE, and here instead by  the ``conserved
quantity" like \eqref{ConQuantity}.

For simplicity, we concentrate on the case $l=2$ here. The proof of
ILP in $H^s$ with negative $s$ and $l\neq 2$ directly follows from
the following argument. In principle, Theorem 1.2 in \cite{FngWng1}
can also be covered by the argument here for the $l\neq 2$ cases and
we'll not exploit it further here.

Note that the ODE part in $t$ for \eqref{focusing-eqn} is
\begin{equation}\label{tode} \pt^2 u=|u|^k |\partial_t
u|\partial_t u\ ,\end{equation} and we have the ``conserved
quantity" for \eqref{tode}(with $u\ge 0$, and $\pt u >
0$)\begin{equation} \ln \pt u-\frac{u^{k+1}}{k+1}.
\end{equation}
If we assign the data $(0,1)$ for \eqref{tode}, then we get a
solution $u_T(t)$ defined on $t\in[0,T)$ with $0<T<\infty$ such
that$$u,\ \pt u,\ \pt^2 u \nearrow \infty\ \ \ \ \rm as\
t\rightarrow T\ .$$ Then for any $a>0$, $u_a(t):=u_T(\frac{T t}{a})$
is the solution of \eqref{tode} with the data $(0,\frac{T}{a})$.
Denote by $\dot{T}_a^s$(or $T^s_a$) the lifespan of the solution of
\eqref{focusing-eqn} with data
$(0,g_a):=(0,\frac{T}{a}\phi(\frac{x}{a}))$ in $\dot{H}^s$(or
$H^s$), where $\phi\in\mathcal{C}_0^\infty$ such that $\phi=1$ on
$B(0,1+d)$ with $d>0.$

If the solution space is $\dot{H}^s$ with $|s|<n/2$, we claim that
\begin{equation} \dot{T}_a^s\leq a.\end{equation} Since, otherwise,
 for $t\in[0,a]$ and $x\in B_{(1+d)a-t}$, the solution $u(t,x)$
equals $u_a(t)$. Note that for $|s|<n/2$, the generalized Leibnitz
rule yields that
\begin{equation}\label{LeibHom}   \|f u\|_{\dot{H}^s}\leq
C\|f\|_{\dot{H}^{n/2}\cap L^\infty}\|u\|_{\dot{H}^s}
\end{equation}
Thus we have
\begin{eqnarray*}\|u(t,x)\|_{\dot{H}_x^s}&\geq& c\|h\|^{-1}_{\dot{H}^{n/2}\cap
L^\infty}\|h(\frac{x}{(1+d)a-t})u(t,x)\|_{\dot{H}_x^s}\\
&=& c u_a(t) ((1+d)a-t)^{\frac{n}{2}-s}\rightarrow\infty
\end{eqnarray*}for $h$ supported in the unit ball $B_1$ as $t\rightarrow a$
from below. On the other hand, for the data, we have(note here that
we choose $\phi$ appropriately as in Section \ref{sltest} such that
$\|\phi\|_{\dot{H}^{s-1}}<\infty$ for any prescribed $s$)
\begin{equation}\|g_a\|_{\dot{H}^{s-1}}=c
a^{\frac{n}{2}-s}\end{equation} Thus by letting $a$ go to zero, we
get the s-ILP of \eqref{focusing-eqn} in $\dot{H}^s$ for $|s|<n/2$.

For the case of $H^s$ with $s<n/2$, by the above result, we have
$$T^s_a\le a$$ for $s\ge 0$, and
$$\|g_a\|_{H^{s-1}}\leq C
(\|g_a\|_{L^2}+\|g_a\|_{\dot{H}^{s-1}})=C(a^{\frac{n}{2}-s}+a^{\frac{n}{2}-1})
$$ for any $s$. For the estimate of lifespan with $s<0$, we
substitute \eqref{LeibHom} by \eqref{Leib} and get the same
estimate. This completes the proof of the $l=2$ case of Theorem
\ref{foc}.

\begin{rem}
Note that  the data we given in this section guarantee the
derivatives of the solution $u_a$ for ODE is nonnegative, one can
change the nonlinear term in \eqref{focusing-eqn} to any reasonable
form such that we also have the ``conserved quantity", say,
$|u|^{k-1} u (\pt u)^2$.
\end{rem}

%\begin{acknowledgements}
%\end{acknowledgements}

% BibTeX users please use
%\bibliographystyle{spmpsci}

\begin{thebibliography}{3}

\bibitem{brenner}
Brenner, P., Kumlin, P.: On wave equations with supercritical
nonlinearities. Arch. Math. {\bf 74}, No.2, 129--147 (2000).

\bibitem{ChCoTa03}  Christ, M., Colliander, J., Tao, T.: Ill-posedness
for nonlinear Schrodinger and wave equations. arxiv:math.AP/0311048

\bibitem{DaGe05}  D'Ancona, P., Georgiev, V.: {Wave maps and ill-posedness
of their Cauchy problem}. In: Reissig, M. (ed.) et al., New trends
in the theory of hyperbolic equations. Basel: Birkh\"{a}user.
Operator Theory: Advances and Applications {\bf 159}, 1--111 (2005)

\bibitem{FngWng1} Fang, D., Wang, C.: Local Well-Posedness and Ill-Posendess
on the Equation of Type $\Box u= u^k (\partial u)^{\alpha}$. Chinese
Ann. Math. Ser. B {\bf 26}, no. 3, 361--378 (2005)

\bibitem{FngWng2} Fang, D., Wang, C.: Sharp Global Existence for Semilinear Wave Equation with Small
 Data. arxiv:math.AP/0612249



\bibitem{lebeau1} Lebeau,  G.: Non linear optic and supercritical wave
equation. Bull. Soc. Roy. Sci. Li\`ege {\bf 70}, no. 4-6, 267--306
(2001)

\bibitem{lebeau2} Lebeau, G.: {Perte de r\'egularit\'e pour les \'equations
d'ondes sur-critiques}, Bull. Soc. Math. France {\bf 133}, no. 1,
145--157 (2005)

\bibitem{93Lind} Lindblad, H.: {A sharp counterexample to the local existence
of low-regularity solutions to nonlinear wave equations}. Duke Math.
J. {\bf 72}, no. 2, 503--539 (1993)

\bibitem{95LdSo} Lindblad, H., Sogge, C.D.:
{On existence and scattering with minimal regularity for semilinear
wave equations}. J. Funct. Anal. {\bf 130}, no. 2, 357--426 (1995)

\bibitem{96Ld} Lindblad, H.:
{Counterexamples to local existence for semi-linear wave equations}.
Amer. J. Math. {\bf 118}, no. 1, 1--16 (1996)


\bibitem{Lind98}Lindblad, H.: {Counterexamples to local existence for
quasilinear wave equations}. Math. Res. Lett. {\bf 5}, no. 5,
605--622 (1998)

\bibitem{95kuksin} Kuksin, S. B.: {On squeezing and flow of energy for
nonlinear wave equations}. Geom. Funct. Anal. {\bf 5}, no. 4,
668--701 (1995)

\end{thebibliography}


\end{document}